\documentclass[12pt]{article}
\usepackage[utf8]{inputenc}
\usepackage{amssymb,amsmath,amsfonts,eucal,mathrsfs,amsthm} 
\usepackage[colorinlistoftodos]{todonotes}
\usepackage[allcolors=blue]{hyperref}
\usepackage{empheq}
\newtheorem{theorem}{Theorem}
\newtheorem{proposition}[theorem]{Proposition}
\newtheorem{lemma}[theorem]{Lemma}

\newtheorem{remark}[theorem]{Remark}

\newtheorem{example}[theorem]{Example}

\theoremstyle{definition}
\newcommand{\R}{\mathbb{R}}

\newcommand{\Sf}{\mathbb{S}}
\newcommand{\C}{\mathbb{C}}

\newcommand{\spa}{\mbox{span}}

\newcommand{\Ima}{\mbox{Im }}

\newcommand{\kerl}{\mbox{ker}}
\newcommand{\grad}{\mbox{grad\,}}

\newcommand{\nab}{\tilde\nabla}

\newcommand{\End}{\mbox{End}}
\newcommand{\Aut}{\mbox{Aut}}

\def\<{{\langle}}
\def\>{{\rangle}}

\def\T{{\cal T}}

\def\n{\nabla}
\def\d{\partial}
\def\a{\alpha}

\def\p{\partial}
\def\be{\begin{equation} }
\def\ee{\end{equation} }
\newcommand{\pu}{\partial_u}
\newcommand{\pv}{\partial_v}
\renewcommand{\gg}{\Gamma^1 }
\newcommand{\gh}{\Gamma^2 }
\def\Ima{Im}
\def\proof{\noindent{\it Proof:  }}
\def\qed{\ifhmode\unskip\nobreak\fi\ifmmode\ifinner
\else\hskip5 pt \fi\fi\hbox{\hskip5 pt \vrule width4 pt
height6 pt  depth1.5 pt \hskip 1pt }}
\makeatletter
\newcommand{\subjclass}[2][]{\let\@oldtitle\@title
\gdef\@title{\@oldtitle\footnotetext{#1 
\emph{Mathematics Subject Classification:} #2}}}
\newcommand{\keywords}[1]{\let\@@oldtitle\@title
\gdef\@title{\@@oldtitle\footnotetext
{\emph{Key words and phrases.} #1.}}}
\makeatother
\begin{document}

\date{}
\title{A construction of Sbrana-Cartan hypersurfaces\\ in the discrete class}
\author{M. Dajczer and M. I. Jimenez}
\keywords{Sbrana-Cartan hypersurfaces, discrete deformations}
\subjclass{53A07, 53B25}
\maketitle

\begin{abstract}  The classical classifications of the locally 
isometrically deformable Euclidean  hypersurfaces obtained by 
U. Sbrana in 1909 and E. Cartan in 1916  includes four classes,
among them the one formed by submanifolds that allow just a single 
deformation. The question of whether these Sbrana-Cartan 
hypersurfaces do, in fact, exist was not addressed by either  
of them. Positive answers to this question were given by 
Dajczer-Florit-Tojeiro in 1998 for the ones called of hyperbolic 
type and by Dajczer-Florit in 2004 when of elliptic type which is 
the other possibility.  In both cases the examples constructed are 
rather special. The main result of this paper yields an abundance 
of examples of hypersurfaces of either type and seems to point 
in the direction of a classification although that goal
remains elusive.

\end{abstract}

An Euclidean hypersurface $f\colon M^n\to\R^{n+1}$, $n\geq 3$, free 
of points where all the sectional curvatures vanish and is locally 
isometrically deformable is called \emph{Sbrana-Cartan}. 
The name is justified because (equivalent) local parametric coordinate 
classifications of this class of submanifolds have been given in terms of 
the nowadays called Gauss parametrization by U. Sbrana \cite{Sb} in 1909 
and by E. Cartan \cite{Ca} in 1916 as envelopes of a two-parameter 
family of affine hyperplanes.  By being isometrically 
deformable we mean that there is at least one other isometric immersion 
of $M^n$ into $\R^{n+1}$ not congruent to $f$ by a rigid motion of the 
ambient space. 
In particular, by the classical Beez-Killing rigidity theorem these 
hypersurfaces must possess two nonzero principal curvatures at any point.

Among the four classes of Sbrana-Cartan hypersurfaces there is the 
one called \emph{surface-like} because the submanifold is a cylinder 
either over a surface in $\R^3$ or over the cone of a surface in 
$\Sf^3\subset\R^4$. Then any isometric deformation of the hypersurface 
is the result of an isometric deformation of the surface. A second class 
consists of the \emph{ruled} hypersurfaces, that is, the ones that carry 
a foliation of codimension one whose leaves are mapped into affine subspaces. 
These submanifolds admit locally as many  isometric deformations as smooth 
real functions on an open interval and all of them carry the 
same rulings. 

Then there are the hypersurfaces that belong to the  important 
additional two classes, namely, the nontrivial ones. One is 
the \emph{continuous class} of hypersurfaces which admit a  
continuous one-parameter family of deformations. For this class, 
the (equivalent) classifications by Sbrana and Cartan are quite 
satisfactory, in particular, because 
they enable the construction of an abundance of explicit examples. 
Finally, there are the hypersurfaces in the remaining 
\emph{discrete class} that admit just  a single deformation. 
The hypersurfaces in this class are of two distinct types called \emph{hyperbolic} and \emph{elliptic} which here are treated 
separately.

The question of whether Sbrana-Cartan hypersurfaces in the discrete 
class do really exist was not at all addressed by either Sbrana or 
Cartan. This is because their parametrically coordinate description  
is quite cumbersome, as seen in detail in the next section. Roughly
speaking, to construct examples one has to search for surfaces for 
which a pair of Christoffel symbols associated to a certain conjugate 
coordinate system has to satisfy a system of second order partial 
differential equations. 

The first positive answer to the above existence question was given 
in \cite{DFT}. 
It was shown that hypersurfaces in the discrete class of hyperbolic type 
can be obtained as the generic transversal intersection of two local 
isometric immersions of $\R^{n+1}$ into $\R^{n+2}$. In fact, the Theorems 
$9$ and $11$ in \cite{DFT} provide a complete description of these 
submanifolds. For an alternative proof of the same result see 
\cite{FF}. Up to now these are the only known examples in this type. 
The existence of examples of hypersurfaces in the discrete class of 
elliptic type was proved in a rather indirect and quite unexpected way. 
In the study of the minimal real Kaehler submanifolds in codimension 
two  in \cite{DF} it was shown that some of these submanifolds are the 
intersection of Sbrana-Cartan hypersurfaces of this type.

In this paper, we construct an abundance of Sbrana-Cartan hypersurfaces 
in the discrete class of both types. Our construction comes out as an 
application of the classification of the infinitesimally bendable 
hypersurfaces, a task initially started by Sbrana \cite{Sb0} and 
then completed by Dajczer-Vlachos in \cite{DV}. For detailed information 
on the subject of infinitesimal bendings
we refer to \cite{DJ1}.
\vspace{1ex}

An \emph{infinitesimal bending} of a hypersurface $f\colon M^n\to\R^{n+1}$
free of flat points is a section $\T$ of $f^*T\R^{n+1}$ that 
satisfies the condition 
\be\label{iif}
\<\T_*X,f_*Y\>+\<f_*X,\T_*Y\>=0
\ee
for any tangent vector fields $X,Y\in\mathfrak{X}(M)$. The bending 
$\T$ is said to be \emph{trivial} if it is the variational vector 
field of a one-parameter variation by isometries of $\R^{n+1}$. 
A hypersurface free of flat points is said to be infinitesimally 
bendable if it admits a nontrivial infinitesimal bending.

The Sbrana-Cartan hypersurfaces that are either surface like, ruled 
or belong the continuous class are infinitesimally bendable whereas 
the ones in the discrete class, as shown below, are not. 
From the classification of the infinitesimally bendable hypersurfaces
in \cite{DV} it follows that the set of the ones that are neither surface 
like or ruled is \emph{much larger} than the set of Sbrana-Cartan 
in the continuous class.

The following is our main result. But first we observe that the 
results in this paper are all of local nature in the sense that 
they hold on connected components of an open dense subset of the 
submanifold.

\begin{theorem}\label{main}
Let $f\colon M^n\to\R^{n+1}$, $n\geq 3$, be an infinitesimally bendable 
hypersurface that is neither surface-like nor ruled on any open subset. 
If $f$ is not a Sbrana-Cartan hypersurface in the continuous 
class then any $f_t=f+t\T$ for $t\in I=(0,\epsilon)$ for some $\epsilon>0$ 
is a  Sbrana-Cartan hypersurface in the discrete class and $f_{-t}$ 
is its isometric deformation. 
\end{theorem}

As seen below Sbrana-Cartan hypersurfaces in the discrete class are not 
infinitesimally bendable which explains why, in the above statement, 
it is only asked the hypersurface not to be in the continuous class.
Moreover, the  conclusion in Theorem \ref{main} does not hold 
if $f$ is allowed to be a Sbrana-Cartan hypersurface in the continuous 
class as shown by Examples \ref{example}. Although we do not know to 
what degree these counter-examples are exceptional we do believe that 
there are no others.

\section{Preliminaries}

In this section, we first give a brief description of the so called 
Gauss parametrization for a class of Euclidean hypersurfaces. For 
plenty of additional information we refer to Chapter $7$ in \cite{DT}. 
After that, we characterize in terms of the Gauss parametrization 
the Sbrana-Cartan hypersurfaces in the continuous class as well as 
the infinitesimally bendable hypersurfaces that are neither surface-like 
nor ruled. In particular, it follows that the later class is much larger 
than the former. For additional details on these subjects, including 
several facts used throughout this paper, we refer to \cite{DFT}, \cite{DT} 
and \cite{DV}.

\subsection{The Gauss parametrization}

Let $f\colon M^n\to\R^{n+1}$, $n\geq 3$, be an isometric immersion 
that at each point has precisely two nonzero principal curvatures. 
The \emph{relative nullity} subspace $\Delta(x)$ of $f$ at $x\in M^n$ 
is the $(n-2)$-dimensional kernel of its second form $A(x)$. 
The relative nullity subspaces form a smooth integrable distribution 
and the totally geodesic leaves of the corresponding 
\emph{relative nullity foliation} are mapped by $f$ into open subsets 
of $(n-2)$-dimensional affine subspaces of $\R^{n+1}$.
\vspace{1ex}

The \emph{Gauss parametrization} of $f$, given in terms of a pair 
$(g,\gamma)$ formed by a surface $g\colon L^2\to\Sf^n$ in the unit 
sphere and a function $\gamma\in C^\infty(L)$, goes as follows: 
Let $U\subset M^n$ be an open saturated subset of leaves of relative 
nullity and $\pi\colon U\to L^2$ the projection onto the 
quotient space. The Gauss map $N$ of $f$ induces an immersion 
$g\colon L^2\to\Sf^n$ in the unit sphere given by $g\circ\pi=N$ 
and the support function $\<f,N\>$, which is constant along the 
relative nullity leaves, induces a function $\gamma\in C^\infty(L)$. 

Let $\Lambda$ denote the normal bundle of $g$ and set $h=i\circ g$ 
where $i\colon\Sf^n\to\R^{n+1}$ is the inclusion. Then the map 
$\psi\colon\Lambda\to\R^{n+1}$ given by
$$
\psi(x,w)=\gamma(x)h(x)+h_*\grad\gamma(x)+i_*w
$$
locally parametrizes (at regular points) the hypersurface $f$ 
in such a way that the fibers of $\Lambda$ are identified with 
the leaves of the relative nullity foliation.

\subsection{The Sbrana-Cartan hypersurfaces}

We describe succinctly in local coordinate terms the Sbrana-Cartan 
hypersurfaces in the discrete and continuous class by means of the 
Gauss parametrization. For proofs and additional information we 
refer to \cite{DFT} and Chapter $11$ of \cite{DT}.
\vspace{2ex}

Let $g\colon L^2\to\Sf^n$ be a surface in the unit sphere.
A system of local coordinates $(u,v)$ on $L^2$ is called
\emph{real-conjugate} for $g$ if its second fundamental form
$\alpha_g\colon TL\times TL\to N_gL$ satisfies
$$
\alpha_g(\d_u,\d_v)=0
$$
where $\pu=\d/\d u$ and $\pv=\d/\d v$. The coordinates are called 
\emph{complex-conjugate} if the condition $\alpha_g(\p_z,\bar\p_z)=0$
holds where  $\p_z=(1/2)(\pu-i\pv)$, that is, if we have 
$$
\alpha_g(\pu,\pu) + \alpha_g(\pv,\pv)=0.
$$

The surface $g\colon L^2\to\Sf^n$ is called \emph{hyperbolic} 
(respectively, \emph{elliptic}) if $L^2$ is endowed with 
real-conjugate (respectively, complex-conjugate) coordinates.
\vspace{1ex}

Set $F=\<\d_u,\d_v\>$ and let $\gg,\gh$ be the Christoffel symbols 
given by $$\n_{\d_u}\d_v=\gg\d_u+\gh\d_v.$$
That $g$ is hyperbolic means the surface 
$h=i\circ g\colon L^2\to\R^{n+1}$ satisfies
\be\label{eq:guv}
h_{uv}-\gg h_u-\gh h_v + Fh=0
\ee
where subscripts indicate partial derivatives.

That the surface $g$ is elliptic means that
$h=i\circ g\colon L^2\to\R^{n+1}$
satisfies
\be\label{eq:guz}
h_{z\bar z}-\Gamma h_z-\bar\Gamma h_{\bar z}+Fh=0,
\ee
where, using the $\C$-linear extensions of the metric of $L^2$ 
and the corresponding connection, the Christoffel symbols are 
given by $\n_{\d_z}\bar\d_z=\Gamma\d_z+\bar\Gamma\bar\d_z$ 
and $F=\<\d_z,\bar\d_z\>$.
\vspace{2ex}

A hyperbolic surface $g\colon L^2\to\Sf^n$ called of 
\emph{first species of real type} if
\be\label{eq:intt} 
\gg_u=\gh_v=2\gg\gh.
\ee
The surface is called  of \emph{second species of real type} 
\index{Surface of \nth{2} species of real type} 
if it is not of first species and the function $\tau=\tau(u,v)$
given by
\be\label{tau}
\tau=\frac{\gg_u-2\gg\gh}{\gh_v-2\gg\gh}
\ee
is positive but not identically one and the necessarily unique 
solution of the system of differential equations
\be\label{eq:siti}
\begin{cases}
\tau_u=2\gh\tau(1-\tau)\vspace{1ex}\\
\tau_v=2\gg(1-\tau).
\end{cases}
\ee 

An elliptic surface $g\colon L^2\to\Sf^n$ is called of 
\emph{first species of complex type} if
\be\label{eq:cos}
\Gamma_z=2\Gamma\bar{\Gamma}(=\bar\Gamma_{\bar z}).
\ee
It is called of \emph{second species of complex type}
\index{Surface of \nth{2} species of complex type} 
if it is not of first species and the function $\rho(z,\bar z)$
with values in the unit sphere and determined by
\be\label{eq:comp}
\Ima(\rho(\Gamma_z-2\Gamma\bar{\Gamma}))=0
\ee
is the necessarily unique non real solution of the differential 
equation
$$
\rho_{\bar z}+\Gamma(\rho-\bar\rho)=0.
$$

The pair $(g,\gamma)$ is called a \emph{hyperbolic pair}
if $g\colon L^2\to\Sf^n$ is an hyperbolic surface  and 
$\gamma\in C^\infty(L)$ satisfies \eqref{eq:guv}, that is, 
the same differential equation than the coordinate functions 
of $h$. Similarly,
the pair $(g,\gamma)$ is called a \emph{elliptic pair} if  
$g\colon L^2\to\Sf^n$ is an elliptic surface and 
$\gamma\in C^\infty(L)$ satisfies \eqref{eq:guz}.

\begin{theorem}\label{Sbrana} {\em (\cite{DFT})}
Let $f\colon M^n\to\R^{n+1}$, $n\geq 3$, be a Sbrana-Cartan 
hypersurface that is neither surface-like nor ruled on any open subset 
of $M^n$. Then, $f$ is parametrized on each connected component 
of an open dense subset of $M^n$ in terms of the Gauss 
parametrization by either an hyperbolic or elliptic pair $(g,\gamma)$, 
where $g\colon L^2\to\Sf^n$ is a surface of first or second species 
of real or complex type. 

Conversely, any simply connected hypersurface parametrized  in terms of 
the Gauss parametrization by such a pair $(g,\gamma)$  is a Sbrana-Cartan 
hypersurface either in the continuous class or in the discrete class, 
according to whether $g$ is, respectively, of first or second species. 
\end{theorem}

\subsection{Infinitesimally bendable hypersurfaces}

We first describe succinctly in local coordinate terms the infinitesimally 
bendable hypersurfaces  by means of the Gauss parametrization. The 
equivalent description in terms of envelopes of hyperplanes is given
in the sequel. For proofs  and additional information we refer to 
\cite{DV} and Chapter $14$ of \cite{DT}.
\vspace{2ex}

An infinitesimal bending $\T$ of $f\colon M^n\to\R^{n+1}$ free of 
flat points is called \emph{trivial} if it is the 
variational vector field of a one-parameter variation by 
isometries of $\R^{n+1}$, namely, if
$$
G(t,x)={\cal C}(t)f(x)+v(t)
$$
where  ${\cal C}\colon I\to O(n+1)$ with $0\in I\subset\R$ 
is a smooth family of orthogonal transformations of $\R^{n+1}$ 
and  $v\colon I\to \R^{n+1}$ a smooth map. Then the variational 
vector field of $G$ at $0\in I$ is  
$$
\T(x)=(\p G/\p t)_{t=0}={\cal D}f(x)+v'(0)
$$
where ${\cal D}={\cal C}'(0)$ is a skew-symmetric linear 
endomorphism of $\R^m$.
\vspace{1ex}

If $\T$ is an infinitesimal bending of $f$ then there is the 
associated smooth 
variation $F\colon\R\times M^n\to\R^{n+1}$  of $f$ by immersions 
$f_t=F(t,\cdot)$ with ${\cal T}$ as the variational vector field 
that is given by $F(t,x)=f(x)+t{\cal T}(x)$. Since
$$
\|f_{t*}X\|^2=\|f_*X\|^2+t^2\|\T_*X\|^2
$$
it is classically said  that the immersions $f_t$
are isometric to $f$ up to the first order.
\vspace{1ex}

The hypersurface  $f\colon M^n\to\R^{n+1}$ is called 
\emph{infinitesimally bendable} if it admits an infinitesimal 
bending that is nontrivial when restricted to any open subset 
of $M^n$. The Sbrana-Cartan 
hypersurfaces that are either surface-like, ruled or belong the 
continuous class are infinitesimally bendable.  In the 
later case, there is the variational vector field of the 
associated one-parameter family of isometric deformations.
\vspace{1ex}

A hyperbolic pair $(g,\gamma)$ is called a \emph{special 
hyperbolic pair} if $g\colon L^2\to\Sf^n$ satisfies
\be\label{integcon}
\gg_u=\gh_v.
\ee
The elliptic pair $(g,\gamma)$ is called a 
\emph{special elliptic pair} if $g\colon L^2\to\Sf^n$ satisfies
\be\label{integcon2}
\Gamma_z=\bar{\Gamma}_{\bar z}, 
\ee
that is, if $\Gamma_z\;\mbox{is real}$.
\vspace{2ex}

In the following result by the infinitesimal bending being unique  
we understand that this is the case up to multiplying 
the bending by a real constant and adding a trivial one.

\begin{theorem}\label{bending} {\em (\cite{DV})} 
Let $f\colon M^n\to\R^{n+1}$, $n\geq 3$, be an infinitesimally
bendable hypersurface that has two nonzero principal curvatures 
everywhere and is neither surface-like nor ruled on any open 
subset of $M^n$. Then $f$ is parametrized on each connected 
component of an open dense  subset of $M^n$ in terms of the 
Gauss parametrization by a special hyperbolic or a special 
elliptic pair.

Conversely, any hypersurface given in terms of the Gauss
parametrization by a special hyperbolic or special elliptic 
pair admits locally a  unique infinitesimal bending.
\end{theorem}

Associated to an infinitesimal bending $\T$ of the hypersurface 
$f\colon M^n\to \R^{n+1}$  with second fundamental form $A$
there is the symmetric tensor $B\in\Gamma(\End(TM))$ 
defined by
$$
\<BX,Y\>=\<(\nab_X\T_*)Y,N\>=\<\nab_X\T_*Y-\T_*\nabla_XY, N\>
$$
which is a Codazzi tensor, that is, $(\nabla_XB)Y=(\nabla_YB)X$ 
and also satisfies the equation
\be\label{eq:awedgeb}
BX\wedge AY-BY\wedge AX=0
\ee
for any $X,Y\in\mathfrak{X}(M)$. That $B$ vanishes says that 
the bending is trivial. Conversely, any nonzero symmetric Codazzi 
tensor $B\in\Gamma(\End(TM))$  satisfying \eqref{eq:awedgeb} 
determines a unique infinitesimal bending of $f$; see Section $2.5$ 
in \cite{DJ1} or Section $5$ in \cite{DV} for details.

We have from \cite{DV} that any infinitesimally bendable hypersurface 
$f\colon M^n\to\R^{n+1}$ can be described as the envelope of a 
two-parameter family of affine hyperplanes. This goes as follows:
Let $U\subset\R^2$ be an open subset endowed with coordinates 
$(u,v)$ and let $\{\varphi_j\}_{0\leq j\leq n+1}$ be a set of 
solutions of the differential equation 
$$
\varphi_{z_1z_2}+M\varphi=0
$$
where $(z_1,z_2)$ can be either $(u,v)$ or $(1/2(u+iv),1/2(u-iv))$ and 
$M\in C^\infty(U)$. Assume that the map 
$\varphi=(\varphi_1,\ldots,\varphi_{n+1})\colon U\to\R^{n+1}$
is an immersion and consider the two-parameter family of affine 
hyperplanes given by
$$
G(u,v)=\varphi_1x_1+\cdots +\varphi_{n+1}x_{n+1}-\varphi_0=0
$$
where $(x_1,\ldots,x_{n+1})$ are the canonical coordinates of 
$\R^{n+1}$. Then $f$  
is the solution of the system of equations $G=G_u=G_v=0$. 
Moreover, the corresponding elliptic pair is given by
$$
g=\frac{1}{\|\varphi\|}(\varphi_1,\ldots,\varphi_{n+1})\;\;\mbox{and}
\;\;\gamma=\frac{\varphi_0}{\|\varphi\|}.
$$

Notice that by this approach we do not require the strong condition that 
a set of solutions of a PDE are the coordinate functions of a surface
in a sphere. We point out that the hypersurface is  Sbrana-Cartan in
the continuous class if the function $\phi=\|\varphi\|^2$
verifies the strong additional condition $\phi_{z_1z_2}=0$;
see \cite{DV} for details.

\begin{proposition}
The Sbrana-Cartan hypersurfaces that belong the discrete class 
are not infinitesimally bendable.
\end{proposition}

\proof For an hyperbolic Sbrana-Cartan hypersurface in the 
discrete class  the function $\tau$ given by \eqref{tau} is 
not identically one.  
In particular, we have that \eqref{integcon} does not hold, 
that is, we have $\gg_u\neq\gh_v$, and thus the hypersurface 
is not infinitesimally bendable.

Similarly, for an elliptic Sbrana-Cartan hypersurface in the discrete 
class  the function $\rho$ given by \eqref{eq:comp} is not real.
In particular, we have that \eqref{integcon2} does not hold,
that is, we have $\Gamma_z\neq\bar{\Gamma}_{\bar z}$,
and again the hypersurface is not infinitesimally bendable.\qed

\begin{remark} {\em In view of the above result we have that an 
infinitesimally bendable hypersurface as in Theorem \ref{bending} 
is not a Sbrana-Cartan hypersurface if either \eqref{eq:intt} in 
the real case or \eqref{eq:cos} in the complex do not hold fully.
}\end{remark}

\section{The main result}

Throughout this section $f\colon M^n\to\R^{n+1}$, $n\geq 3$, 
stands for an infinitesimally bendable hypersurface with infinitesimal 
bending $\T$ that is free of flat points and is neither surface-like 
nor ruled restricted to any open subset of $M^n$.
\vspace{1ex}

We recall that the \emph{splitting tensor} 
$C\colon\Gamma(\Delta)\to\Gamma(\End(\Delta^\perp))$ of
the $(n-2)$-dimensional relative nullity distribution $\Delta$ 
of $f$ is defined by 
$$
C(T,X)=-(\n_XT)_{\Delta^\perp}
$$
and $C_T$ denotes the element of $\Gamma(\End(\Delta^\perp))$ given 
by $C_TX=C(T,X)$.
\vspace{1ex}

The assumption on its geometry gives that $f$ restricted to each 
connected component of an open dense subset of $M^n$ satisfies that 
$C_T$ for some $T\in\Gamma(\Delta)$ has either two nonzero distinct 
real or complex conjugate eigenvalues at any point. In fact, we
have from Proposition $7.4$ in \cite{DT} that $f$ is surface-like 
with respect to the decomposition $TM=\Delta\oplus\Delta^\perp$ 
if the splitting tensor satisfies $C_T\in\spa\{I\}$ for any 
$T\in\Gamma(\Delta)$. If there is  $T\in\Gamma(\Delta)$ such that 
$C_T$ has a single real eigenvalue of multiplicity $2$ then the 
argument in pp. $372-373$ in \cite{DFT} gives that $f$ is ruled. 
\vspace{1ex}

The hypersurface $f\colon M^n\to\R^{n+1}$ is called \emph{hyperbolic} 
(respectively, \emph{elliptic}) if there exists a tensor  
$J\in\Gamma(\Aut(\Delta^\perp))$ that satisfies:
\begin{itemize}
\item[(i)]  $J^2=I$ with $J\neq I$ (respectively, $J^2=-I$), 
\item[(ii)] $\nabla_TJ=0$ for all $T\in\Gamma(\Delta)$, 
\item[\hypertarget{(iii)}{(iii)}] $C_T\in\spa\{I,J\}$ for all 
$T\in\Gamma(\Delta)$ but $C(\Gamma(\Delta))\not\subset\spa\{I\}$.
\end{itemize}
Under the assumptions required for $f$ it follows from Proposition $20$ 
in \cite{DV} that its restriction to any  connected component of an
open dense subset of $M^n$ is either hyperbolic or elliptic. 

\begin{lemma}\label{relnul}
The relative nullity distribution $\Delta$ of $f$ is contained at any 
point of $M^n$ in the relative nullity $\Delta_t$ of $f_t=f+t\T$ for any $t\in\R$.
\end{lemma}

\proof From the classical Beez-Killing rigidity theorem and either Theorem 
$1$ in \cite{DR} or Proposition $14.3$ in \cite{DT} we have that 
the second fundamental 
form $A_t$ of $f_t$ has rank at most $2$ for any $t\in\R$. 
The unit normal vector field $N_t$ of $f_t$ decomposes as
$$
N_t=f_*Z_t+b N
$$
where $Z_t\in \mathfrak{X}(M)$ and $b=b(t,x)=\<N_t,N\>$.
Let $\mathcal{Y}\in\mathfrak{X}(M)$ be given by
$$
\<\mathcal{Y},X\>+\<N,\T_*X\>=0
$$
for all $X\in\mathfrak{X}(M)$. Then using \eqref{iif} we have
\begin{align*}
0=\<N_t,f_{t*}X\>&=\<f_*Z_t+bN,f_*X+t\T_*X\>\\
&=\<Z_t,X\>+t\<f_*Z_t,\T_*X\>+tb\<N,\T_*X\>\\
&=\<f_*Z_t-t\T_*Z_t-tbf_*\mathcal{Y},f_*X\>,
\end{align*}
that is,
\be\label{z1}
(f_*Z_t-t\T_*Z_t-tbf_*\mathcal{Y})_{f_*TM}=0.
\ee

Let $L_0\in\Gamma(\End(TM))$ be defined by
$$
f_*L_0X=\pi\T_*X
$$
where $\pi\colon \R^{n+1}\to f_*TM$ is the orthogonal projection. 
Then \eqref{z1} reads as
\be\label{reads}
(I-tL_0)Z_t=tb\mathcal{Y}.
\ee
Since $L_0$ is skew-symmetric by \eqref{iif} then 
\kerl($I-tL_0)=0$. If $S_t=(I-tL_0)^{-1}$ then 
\be\label{Zt}
Z_t=tbS_t\mathcal{Y}
\ee
where $b\neq 0$ since otherwise $S_t$ would have 
a nontrivial kernel.

Given $T\in\Gamma(\Delta)$ and being $\pi$ parallel along the 
leaves of relative nullity, we have
$$
(\nab_T f_*L_0)X=\nab_Tf_*L_0X-f_*L_0\nabla_TX=\pi(\nab_T\T_*)X.
$$
We obtain from \eqref{eq:awedgeb} that $\Delta\subset\ker B$.
On the other hand, we have by Proposition $2.2$ in \cite{DJ1} that 
$$
(\nab_T\T_*)X=\<AT,X\>f_*\mathcal{Y}+\<BT,X\>N=0.
$$
Therefore
\be\label{derl0}
\nab_T f_*L_0=0
\ee
for any $T\in\Gamma(\Delta)$.
From $(2.16)$ in \cite{DJ1}  or $(13)$ in \cite{DJ} we obtain
$$
\nab_X\mathcal{Y}=-f_*BX-\T_*AX
$$
and hence
\be\label{dery}
\nab_T\mathcal{Y}=0
\ee
for any $T\in\Gamma(\Delta)$.
Taking the covariant derivative of \eqref{reads} with respect to
$T$ and using \eqref{derl0} and \eqref{dery} it follows that
$$
\nab_Tf_*Z_t=tT(b)S_t\mathcal{Y}.
$$
Then
$$
\nab_TN_t=T(b)(tS_t\mathcal{Y}+N).
$$
Since $b\neq0$ we obtain from \eqref{Zt} that $\nab_TN_t$ is a 
multiple of $N_t$ which is of unit length. Hence $\nab_TN_t=0$ for  
any $T\in\Gamma(\Delta)$, and this proves that $\Delta\subset\Delta_t$
for any $t\in\R$.
\vspace{2ex}\qed

In the sequel, we assume that $\epsilon>0$ is small enough so that 
$\kerl A(t)=\Delta_t=\Delta$ for any $t\in (-\epsilon,\epsilon)$. 
Notice that the orthogonal complement $\Delta_t^\perp$ does not have 
to coincide with $\Delta^\perp$. Let $C^t$ denote the splitting tensor 
of $\Delta$ with respect to the metric determined by $f_t$. 
An argument of continuity  applied to $C^t$ gives that $f_t$ 
for some $\epsilon>0$ is neither surface-like nor ruled on an 
open subset of $M^n$. 

\begin{lemma} \label{remain}
If $f$ is either hyperbolic or elliptic then the immersions $f_t$, 
$t\in [0,\epsilon)$ for small enough $\epsilon>0$ remain hyperbolic 
or elliptic with respect to a tensor $J_t\in\Gamma(\End(\Delta_t^\perp))$.
\end{lemma}

\proof Since the $f_t$ are Sbrana-Cartan hypersurfaces that for $t$ 
small enough are neither surface-like nor ruled, then from the proof of 
Theorem $11.16$ in\cite{DT}  we have that they are either elliptic 
or hyperbolic with respect to $J_t\in\Gamma(\End(\Delta_t^\perp))$ 
for each $t$. From that proof we also have that the splitting tensor 
$C^t$ for each $t$ satisfies $C^t_T\in\spa\{I,J_t\}$ for all 
$T\in\Gamma(\Delta)$. Then, by continuity, the conditions 
$J_t^2=I$ or $J_t^2=-I$ remain for small values of $t$.
\vspace{2ex}\qed
 
The Gauss map $N_t$ of $f_t$ determines an immersion 
$g_t\colon L_t^2\to\Sf^n$ on the quotient space of leaves of 
relative nullity $L_t^2$ into the unit sphere. 
By Lemma~\ref{relnul} we have  that $L_t=L$ and hence the 
family of immersions $g_t\colon L^2\to\Sf^n$ depends smoothly 
on the parameter $t$.
Since the tensors $J_t\in\Gamma(\End(\Delta_t^\perp))$ have to
satisfy the property \hyperlink{(iii)}{$(iii)$} it follows that 
the dependence on the parameter $t$ is smooth. 

By Proposition $11.11$ in \cite{DT} we have that the tensor 
$J_t$ is the horizontal lifting of a tensor 
$\bar{J}_t\in\Gamma(\End(TL))$,
i.e., $\bar{J}_t\circ\pi_*=\pi_*\circ J_t$, such that
$\bar{J}_t^2=I$ ($\bar{J}_t\neq I$) or $\bar{J}_t^2=-I$ 
according to $J_t$. 
Moreover, the immersion $g_t$ is hyperbolic or elliptic with
respect to $\bar{J}_t$, accordingly, which in this case means 
that the second fundamental form of $g_t$ satisfies
$$
\a^{g_t}(\bar{J}_tX,Y)=\a^{g_t}(X,\bar{J}_tY)
$$
for any $X,Y\in\mathfrak{X}(L)$.

Finally,  from the proof of Theorem $11.16$ in \cite{DT} we have
that a deformation of $f_t$ for $t\neq 0$ is determined 
by a tensor tensor $\bar{D}_t\in\Gamma(\End(TL))$ satisfying:
\begin{itemize}
\item[(i)] $\bar{D}_t\in\spa\{\bar{I},\bar{J}_t\}$, 
$\bar{D}_t\neq\pm\bar{I}$,
\item[(ii)] $\det\bar{D}_t=1$,
\item[(iii)] $\bar{D}_t$ is a Codazzi tensor with respect to 
the metric induced by $g_t$.
\end{itemize}

\subsection{The hyperbolic case}

We focus on the case when $f$ is hyperbolic, that is, 
we have $\bar{J}_t^2=\bar{I}$ for $t\in I=[0,\epsilon)$.

\begin{lemma}\label{frame1}
There is locally a smooth one-parameter family of tangent 
frames $\{U^t,V^t\}$ by vectors  fields $U^t,V^t\in\mathfrak{X}(L)$ 
of unit length with respect to the metric of $g_t$ that 
are eigenvectors of $\bar{J}_t$ associated to the eigenvalues 
$1$ and $-1$, respectively.
\end{lemma}

\proof Let $U,V\in\mathfrak{X}(L)$ be a frame of eigenvectors of 
$\bar{J}$ of unit length with respect to the metric induced by 
$g=g_0$ associated to the eigenvalues $1$ and $-1$, respectively. 
With respect to this frame we have
$$
\bar{J}_t=
\begin{bmatrix}
a & c\\
b & d
\end{bmatrix}.
$$
That $\bar{J}_t^2=\bar{I}$ is equivalent to $a+d=0$ and 
$1-a^2=bc$. If $\mu,\nu\in C^\infty(I\times L)$ are given by 
$\mu=b/(1-d)$ and $\nu=-c/(1+a)$ then the vector fields 
$U+\mu V$ and $V+\nu U$ are eigenvectors of $\bar{J}_t$. 
We obtain the desired frames by normalizing them.
\vspace{2ex}\qed

Let $U^t,V^t$ be the frame given by Lemma \ref{frame1}. Since 
$\det\bar{D}_t=1$ there is $\theta\in C^\infty(I\times L)$ such 
that $\bar{D}_t$ in terms of $\{U^t,V^t\}$ is of the form
$$
\bar{D}_t=
\begin{bmatrix}
\theta & 0\\
0 & \theta^{-1}
\end{bmatrix}.
$$

The following calculations are all done for a fixed $t$ 
that is omitted in the sequel for simplicity of notation. 
Since $U,V$ are unit vector fields we have  
$\<\nabla_UV,V\>=0=\<\nabla_VU,U\>$. Thus 
\be\label{form}
\nabla_UV=\Lambda_1U-\Lambda_1FV\;\;\mbox{and}\;\;
\nabla_VU=-\Lambda_2FU+\Lambda_2V
\ee
where $F=\<U,V\>$. Then the Codazzi equation for $\bar D$ is 
equivalent to the system of equations
$$
\begin{cases}
U(\theta^{-1})=\Lambda_2(\theta-\theta^{-1})\vspace{1ex}\\
V(\theta)=\Lambda_1(\theta^{-1}-\theta).
\end{cases}
$$
Multiplying the first equation by $2\theta^3$, the 
second by $2\theta$  and setting $\tau=\theta^2$ yields
$$
\begin{cases}
U(\tau)=2\Lambda_2\tau(1-\tau)\vspace{1ex}\\
V(\tau)=2\Lambda_1(1-\tau).
\end{cases}
$$
Then any positive solutions of this system other than $\tau=1$ 
gives a Codazzi tensor $\bar{D}$ with $\det\bar{D}=1$. The 
integrability condition of the system is 
$$
U(\Lambda_1)-\Lambda_1\Lambda_2+F(\Lambda_1)^2
-\tau(V(\Lambda_2)-\Lambda_1\Lambda_2+F(\Lambda_2)^2)=0
$$
that is trivially satisfied if 
\be\label{fsrt1}
\Lambda_1\Lambda_2=U(\Lambda_1)+F(\Lambda_1)^2=
V(\Lambda_2)+F(\Lambda_2)^2.
\ee

\begin{lemma}
The hyperbolic surface $g_t\colon L^2\to\Sf^n$ is of the first 
species of real type if and only if the condition \eqref{fsrt1} holds.
\end{lemma}

\proof  Let $(u,v)$ be a coordinate system such that 
$\partial_u=aU$ and $\partial_v=bV$ for  $a,b\in C^{\infty}(L)$.
Then
$$
\begin{cases}
\nabla_{\partial_u}\partial_v=b \Lambda_1\partial_u
+a(U(\log b)-F\Lambda_1)\partial_v\vspace{1ex}\\
\nabla_{\partial_v}\partial_u=b(V(\log a)
-F\Lambda_2)\partial_u+a\Lambda_2\partial_v.
\end{cases}
$$
Since $[\partial_u,\partial_v]=0$ we have
\be\label{derab}
\Lambda_1=V(\log a)-F\Lambda_2\;\;\mbox{and}
\;\; \Lambda_2=U(\log b)-F\Lambda_1,
\ee
and hence
$$
\nabla_{\partial_u}\partial_v
=b\Lambda_1\partial_u+a\Lambda_2\partial_v.
$$
Hence the surface is of the first species of real type if 
$$
\partial_u(b\Lambda_1)=\partial_v(a\Lambda_2)
=2ab\Lambda_1\Lambda_2
$$
which using \eqref{derab} is verified to be equivalent 
to \eqref{fsrt1}.
\vspace{2ex}\qed

Next we prove Theorem \ref{main} in the case when $f$ is hyperbolic.
\vspace{2ex}

\proof We assume that the hypersurface is parametrized by a special 
hyperbolic pair. Let $g_t\colon L^2_t\to\Sf^n$ be the one-parameter 
family of surfaces determined by the Gauss maps of $f_t$ and 
$t\in [0,\epsilon)$ for $\epsilon>0$ given by Lemmas \ref{relnul}
and \ref{remain}. For $U^t$ and $V^t$ as in Lemma \ref{frame1} 
we have from \eqref{form} that $\Lambda_1$ and $\Lambda_2$ 
are smooth on the parameter $t$. Hence, if there is a sequence 
$\{t_n\}\to 0$ in $[0,\epsilon)$ such that for each $t_n$ the 
surface $g_{t_n}$ is of the first species of real type then 
\eqref{fsrt1} holds for each $t_n$. By continuity \eqref{fsrt1}
also holds for $t=0$ in contradiction with the assumption on $f$ 
not to be Sbrana-Cartan.\qed

\subsection{The elliptic case}

Next we treat the case when $f$ is elliptic, that is, when 
$\bar{J}^2_t=-\bar{I}$ for $t\in I$. 
\vspace{1ex}

Let $\{U^t,V^t=\bar{J}_tU^t\}$ be a smooth 
one-parameter family of tangent frames. Since
$\bar{D}_t\in\spa\{\bar{I},\bar{J}\}$ and $\det\bar{D}_t=1$ then
$\bar{D}_t$ in that frame is of the form
$$
\bar{D}_t=\begin{bmatrix}
\cos\theta&-\sin\theta\\
\sin\theta&\cos\theta
\end{bmatrix}
$$
where $\theta\in C^\infty(L)$. We consider the complexified 
tangent bundle $TL^\C$ as well as the $\C$-linear extension of the 
tensors $\bar{D}_t$ and $\bar{J}_t$ denoted  
equally. Then the vector field $Z^t=1/2(U^t-iV^t)$ and $\bar{Z^t}$ 
are eigenvectors of $\bar{D}_t$ with respect to the eigenvalues 
$\rho=\cos\theta+i\sin\theta$ and $\bar{\rho}$, respectively. 

The following calculations are all done for a fixed $t$ that is
omitted in the sequel for simplicity of notation. Now 
we consider the $\C$-bilinear extensions of the metric and 
the Levi-Civita connection. Then we write
\be\label{chrissym2}
\nabla_Z\bar{Z}=\Lambda_1Z+\Lambda_2\bar{Z}
\ee
and observe that $\nabla_{\bar{Z}}Z=\overline{\nabla_Z\bar{Z}}$.

The Codazzi equation for $\bar{D}$ evaluated in $Z$ and $\bar{Z}$ 
is equivalent to 
\be\label{Codazzi2}
Z(\bar{\rho})=\bar{\Lambda}_1(\rho-\bar{\rho}).
\ee
Since $\rho\bar{\rho}=1$, then that $Z(\rho\bar\rho)=0$ together 
with \eqref{Codazzi2} yields
$$
Z(\rho)=\rho^2\bar{\Lambda}_1(\bar{\rho}-\rho)=-\rho^2Z(\bar{\rho}).
$$
Now a straightforward computation of the integrability 
condition for \eqref{Codazzi2} gives
$$
Z(\Lambda_1)-\Lambda_1\bar{\Lambda}_1-\Lambda_1\Lambda_2
=\rho^2(\bar{Z}(\bar{\Lambda}_1)-\Lambda_1\bar{\Lambda}_1
-\bar{\Lambda}_1\bar{\Lambda}_2).
$$
Then this condition holds trivially if 
\be\label{fsct1}
Z(\Lambda_1)=\Lambda_1\bar{\Lambda}_1+\Lambda_1\Lambda_2.
\ee

\begin{lemma}
The elliptic surface $g_t\colon L^2\to\Sf^n$ is of the first species 
of complex type if and only if the condition \eqref{fsct1} holds.
\end{lemma}

\proof   Let $(u,v)$ be complex conjugate coordinates such that
$\frac{1}{2}(\partial_u-i\partial_v)=\partial_z=\phi(z)Z$ and
write
$$
\nabla_{\partial_z}\partial_{\bar{z}}
=\Gamma\partial_z+\bar{\Gamma}\partial_{\bar{z}}.
$$
A straightforward computation, using \eqref{chrissym2} gives
$$
\begin{cases}
 \nabla_{\partial_z}\partial_{\bar{z}}=\phi\bar{\phi}\Lambda_1Z
+(\phi Z(\bar{\phi})+\phi\bar{\phi}\Lambda_2)\bar{Z}\vspace{1ex}\\
\nabla_{\partial_{\bar{z}}}\partial_z=(\bar{\phi}\bar{Z}(\phi)
+\phi\bar{\phi}\bar{\Lambda}_2)Z
+\phi\bar{\phi}\bar{\Lambda}_1\bar{Z}.
\end{cases}
$$
Hence $\nabla_{\partial_z}\partial_{\bar{z}}
=\nabla_{\partial_{\bar{z}}}\partial_z$ yields
\be\label{derphi}
Z(\bar{\phi})+\bar{\phi}\Lambda_2=\bar{\phi}\bar{\Lambda}_1
\ee
and 
\be\label{delta}
\Gamma=\bar{\phi}\Lambda_1.
\ee
Recall that the surface is of the first species of complex 
type if
\be\label{fsct2}
\Gamma_z=2\Gamma\bar{\Gamma}.
\ee
Since \eqref{delta} gives
$$
\Gamma_z=\phi Z(\bar{\phi})\Lambda_1
+\phi\bar{\phi}Z(\Lambda_1)
$$
we have by \eqref{derphi} that \eqref{fsct1} is equivalent 
to \eqref{fsct2}.\qed

\begin{lemma}\label{frame2}
Assume that  $g_0=g$ is not a minimal surface. Then there is 
a local smooth one-parameter family of vector fields  
$U^t,V^t\in\mathfrak{X}(L)$ of unit norm with respect to the metric 
induced by $g_t$ such that $\bar{J}^tU^t=V^t$, $\bar{J}^tV^t=-U^t$
for $t\in (-\epsilon,\epsilon)$ and small $\epsilon>0$.
\end{lemma}

\proof Take a local frame $U^t,V^t=\bar{J}^tU^t\in\mathfrak{X}(L)$. 
Again we omit for simplicity of notation the index $t$. 
We search for functions $d,e\in C^\infty(I\times L)$ such that
$\tilde{U}=dU+eV$ and $\tilde{V}=\bar{J}\tilde{U} =dV-eU$ have
unit length, that is, we look for solutions of the system
$$
\begin{cases}
d^2\|U\|^2+2de\<U,V\>+e^2\|V\|^2=1\vspace{1ex}\\
e^2\|U\|^2-2de\<U,V\>+d^2\|V\|^2=1.
\end{cases}
$$
Equivalently, we want solutions of the equations
$$
(d^2+e^2)(\|U\|^2+\|V\|^2)=2
$$
and 
$$
(d^2-e^2)\|U\|^2+4de\<U,V\>-(d^2-e^2)\|V\|^2=0.
$$
The later is the equation of second degree in $d/e$ given by
$$
\frac{d^2}{e^2}(\|U\|^2-\|V\|^2)
+4\dfrac{d}{e}\<U,V\>+\|V\|^2-\|U\|^2=0
$$
which is not trivial since  the surfaces $g_t$ are not 
minimal for $t$ small.\vspace{2ex}\qed

Finally, we prove Theorem \ref{main} in the case when $f$ is  
elliptic.
\vspace{2ex}

\proof  We assume that the hypersurface is parametrized by a special 
elliptic pair and argue similarly to the hyperbolic case but now 
we use Lemmas \ref{frame2} and \eqref{fsct1}. Notice that since $g$ 
is not of the first species of complex type then it is not a minimal 
surface.\qed

\begin{example}\label{example}
{\em Any simply connected isometric minimal immersion 
$f\colon M^{2n}\to\R^{2n+1}$, $n\geq2$, of a nonflat Kaehler manifold 
is a Sbrana-Cartan hypersurface in the continuous class since there is 
an associated one-parameter family $f^{\theta}\colon M^{2n}\to\R^{2n+1}$, 
$\theta\in [0,\pi)$ of isometric submanifolds given by 
$f^{\theta}=\cos\theta f+\sin\theta\bar{f}$ all having the same Gauss 
map; see Theorem $15.8$ in \cite{DT}. Moreover, the conjugate submanifold 
$\bar f=f^{\pi/2}$ satisfies $\<f_*X,\bar{f}_*X\>=0$ for any 
$X\in\mathfrak{X}(M)$, and hence is an infinitesimal bending of $f$;  
see Proposition $15.9$ in \cite{DT}.  Then the submanifold 
$f_t=f+t\bar{f}$ for any $t\in\R$ is also a Sbrana-Cartan hypersurface 
in the continuous class since it is homotetic to an element in the 
associated family to $f$.
}\end{example}

\section*{Acknowledgment}

The first author thanks the Mathematics Department of the University
of Murcia where most of this work was developed for the kind 
hospitality during his visit.
The second author is supported by CAPES-PNPD Grant 88887.469213/2019-00.

\medskip

This research is part of the grant PID2021-124157NB-I00, funded by\\
MCIN/ AEI/10.13039/501100011033/ ``ERDF A way of making Europe".

\noindent Marcos Dajczer\\
IMPA -- Estrada Dona Castorina, 110\\
22460--320, Rio de Janeiro -- Brazil\\
e-mail: marcos@impa.br

\bigskip

\noindent Miguel Ibieta Jimenez\\
Universidade de S\~ao Paulo\\
Instituto de Ci\^encias Matem\'aticas e de Computa\c c\~ao\\
Av. Trabalhador S\~ao Carlense 400\\
13566--590, S\~ao Carlos -- Brazil\\
e-mail: mibieta@impa.br
\end{document}